\documentclass[12pt]{article}
\usepackage{amsfonts}
\usepackage{amsmath}
\usepackage{graphicx}

\input tcilatex
\input tcilatex
\begin{document}

\author{I. Dag$^{1}$, O Ersoy$^{1}$ and \"{O} Ka\c{c}maz$^{2}$ \\
%EndAName
Eski\c{s}ehir Osmangazi University, \\
Faculty of Science and Art,\\
Department of Mathematics-Computer, Eski\c{s}ehir, Turkey$^{1}$,\\
Open Education Faculty, Distance Education Division. \\
Anadolu University Turkey$^{2}$, }
\title{\textbf{The Trigonometric Cubic B-spline Algorithm for Burgers'
Equation}}
\maketitle

\begin{abstract}
The cubic Trigonometric B-spline(CTB) functions are used to set up the
collocation method for finding solutions of the Burgers' equation. The
effect of the CTB in the collocation method is sought by studying two text
problems. The Burgers' equation is fully-discretized using the
Crank-Nicholson method for the time discretizion and CTB function for
discretizion of spatial variable. Numerical examples are performed to show
the convenience of the method for solutions of Burgers equation

\end{abstract}

\section{Introduction}

Since the introduction of the Burger's equation by Bateman\cite{bat}, many
authors have used variety of numerical methods in attempting to solve the
Burger's equation.\ Various forms of the finite element methods are
constructed to compute the Burger's equation numerically\cite%
{ad,ad1,nu2,sku,ad5,ofem3,BSpGal,qq,zi}. The spline collocation procedures
are also presented for getting solutions of the Burger's equation\cite%
{rt,er,ras,CBSpCol,irk,ram,QBSpCol2,QBSpCol,dag,QuinCol}. The spline
functions are wished to be accompanied to the numerical method to solve the
differential equations since the resulting matrix system is always diagonal
and can be solved easily and approximate solutions having the accuracy of
the degree less than the degree of the spline functions can be set up. High
order continuous differentiable approximate solutions can be produced for
the differential equations of higher order. The numerical procedure for
nonlinear evolution equations based on the the B-spline collocation method
have been increasingly applied to various fields of science. However
application of the CTB collocation method to non linear evolution problems
are a few in comparison with the method the collocation based on the
B-spline functions.

The numerical methods for solving a type of ordinary differential equations
with quadratic and cubic CTB are given by A. Nikolis in the papers\cite%
{g1,g2}. The linear two-point boundary value problems of order two are
solved using cubic CTB interpolation method \cite{we}. The another numerical
method employed the cubic CTB are set up to solve a class of linear
two-point singular boundary value problems in the study\cite{w1}. Very
recently a collocation finite difference scheme based on new cubic CTB is
developed for the numerical solution of a one-dimensional hyperbolic
equation (wave equation) with non-local conservation condition\cite{aa}. \ A
new two-time level implicit technique based on the cubic CTB is proposed for
the approximate solution of the nonclassical diffusion problem with nonlocal
boundary condition in the study\cite{ma}. \ Some researches have established
types of the B-spline finite element approaches for solving the Burger's
equation \ but \ not with CTB as far as we know the literature.

In this paper, CTB are used to establish a collocation method \ and then the
suggested numerical method is applied to find the numerical solutions of
Burger's equation. It is also well-known that this problem arise in many
branch of the science and the development of the numerical methods for the
Burger's equation have been attracted for finding the steep front solutions.
The use of the lower \ viscosity in the Burger's equation cause the
appearance of the steep front and shock wave solutions. That makes
difficulty in modelling solutions of the Burger's equation when solved
numerically. So many authors have developed various kinds of numerical
scheme in computing the equation effectively for small values of the
viscosity.

We consider the Burger's equation

\begin{equation}
U_{t}+UU_{x}-\lambda U_{xx}=0,\text{ }a\leq x\leq b,\text{ }t\geq 0
\label{1}
\end{equation}%
with appropriate initial conditions and the boundary conditions:$%
U(x,0)=f(x), $ $a\leq x\leq b,U(a,t)=U_{a}$, $U(b,t)=U_{b}$ where subscripts 
$x$ and $t$ denote differentiation, $\lambda =\dfrac{1}{\func{Re}}>0$ and $%
\func{Re}$ is the Reynolds number characterizing the strength of viscosity. $%
U_{a},$ $U_{b} $ are the constants and $u=u(x,t)$ is a sufficiently
differentiable unknown function and $f(x)$ is a bounded function.

The implementation of the proposed scheme is given in the second section.
Two classical text problems are dealt with to show the robustness of the
scheme.

\section{Cubic Trigonometric B-spline Collocation Method}

Consider a uniform partition of the problem domain $[a=x_{0},b=x_{N}]$ at
the knots $x_{i},i=0,...,N$ with mesh spacing $h=(b-a)/N.$ On this partition
together with additional knots $x_{-1},x_{0},x_{N+1},x_{N+2},x_{N+3}$
outside the problem domain, $CTB_{i}${}$(x)$ can be defined as

\begin{equation}
CTB_{i}(x)=\frac{1}{\theta }\left \{ 
\begin{tabular}{ll}
$\omega ^{3}(x_{i-2}),$ & $x\in \left[ x_{i-2},x_{i-1}\right] $ \\ 
$\omega (x_{i-2})(\omega (x_{i-2})\phi (x_{i})+\phi (x_{i+1})\omega
(x_{i-1}))+\phi (x_{i+2})\omega ^{2}(x_{i-1}),$ & $x\in \left[ x_{i-1},x_{i}%
\right] $ \\ 
$\omega (x_{i-2})\phi ^{2}(x_{i+1})+\phi (x_{i+2})(\omega (x_{i-1})\phi
(x_{i+1})+\phi (x_{i+2})\omega (x_{i})),$ & $x\in \left[ x_{i},x_{i+1}\right]
$ \\ 
$\phi ^{3}(x_{i+2}),$ & $x\in \left[ x_{i+1},x_{i+2}\right] $ \\ 
\multicolumn{1}{c}{} & \multicolumn{1}{c}{$0,\text{otherwise}$}%
\end{tabular}%
\right.  \label{r2}
\end{equation}%
where $\omega (x_{i})=\sin (\frac{x-x_{i}}{2}),\phi (x_{i})=\sin (\frac{%
x_{i}-x}{2}),\theta =\sin (\frac{h}{2})\sin (h)\sin (\frac{3h}{2}).$

$CTB_{i}(x)$ are twice continuously differentiable piecewise trigonometric
cubic B-spline on the interval$[a,b]$. \ The iterative formula \bigskip

\begin{equation}
T_{i}^{k}(x)=\frac{\sin (\frac{x-x_{i}}{2})}{\sin (\frac{x_{i+k-1}-x_{i}}{2})%
}T_{i}^{k-1}(x)+\frac{\sin (\frac{x_{i+k}-x}{2})}{\sin (\frac{x_{i+k}-x_{i+1}%
}{2})}T_{i+1}^{k-1}(x),\text{ }k=2,3,4,...  \label{r3}
\end{equation}%
gives the cubic B-spline trigonometric functions starting with the
CTB-splines of order $1$

\begin{equation*}
T_{i}^{1}(x)=\left \{ 
\begin{tabular}{c}
$1,\ x\in \lbrack x_{i},x_{i+1})$ \\ 
$0$ \  \ ,otherwise.%
\end{tabular}%
\right.
\end{equation*}

Each $CTB_{i}(x)$ is twice continuously differentiable and the values of $%
CTB_{i}(x),CTB_{i}^{^{\prime }}(x)$ and $CTB_{i}^{^{\prime \prime }}(x)$ at
the knots $x_{i}$ 's can be computed from Eq.(\ref{r3}) as

\begin{equation*}
\begin{tabular}{l}
Table 1: Values of $B_{i}(x)$ and its principle two \\ 
derivatives at the knot points \\ 
\begin{tabular}{|l|l|l|l|}
\hline
& $T_{i}(x_{k})$ & $T_{i}^{\prime }(x_{k})$ & $T_{i}^{\prime \prime }(x_{k})$
\\ \hline
$x_{i-2}$ & 0 & 0 & 0 \\ \hline
$x_{i-1}$ & $\sin ^{2}(\frac{h}{2})\csc \left( h\right) \csc (\frac{3h}{2})$
& $\frac{3}{4}\csc (\frac{3h}{2})$ & $\frac{3(1+3\cos (h))\csc ^{2}(\frac{h}{%
2})}{16\left[ 2\cos (\frac{h}{2})+\cos (\frac{3h}{2})\right] }$ \\ \hline
$x_{i}$ & $\frac{2}{1+2\cos (h)}$ & 0 & $\frac{-3\cot ^{2}(\frac{3h}{2})}{%
2+4\cos (h)}$ \\ \hline
$x_{i+1}$ & $\sin ^{2}(\frac{h}{2})\csc \left( h\right) \csc (\frac{3h}{2})$
& -$\frac{3}{4}\csc (\frac{3h}{2})$ & $\frac{3(1+3\cos (h))\csc ^{2}(\frac{h%
}{2})}{16\left[ 2\cos (\frac{h}{2})+\cos (\frac{3h}{2})\right] }$ \\ \hline
$x_{i+2}$ & 0 & 0 & 0 \\ \hline
\end{tabular}%
\end{tabular}%
\end{equation*}%
$CTB_{i}(x)$ , $i=-1,...,N+1$ are a basis for the trigonometric spline
space. An approximate solution $U_{N}$ to the unknown $U$ is written in
terms of the expansion of the CTB as

\begin{equation}
U_{N}(x,t)=\sum_{i=-1}^{N+1}\delta _{i}CTB_{i}(x)  \label{r4}
\end{equation}%
where $\delta _{i}$ are time dependent parameters to be determined from the
collocation points $x_{i},i=0,...,N$ and the boundary and initial
conditions. The nodal values $U$ and its first and second derivatives at the
knots can be found from the (\ref{r4}) as 
\begin{equation}
\begin{tabular}{l}
$U_{i}=\alpha _{1}\delta _{i-1}+\alpha _{2}\delta _{i}+\alpha _{1}\delta
_{i+1}$ \\ 
$U_{i}^{\prime }=\beta _{1}\delta _{i-1}+\beta _{2}\delta _{i+1}$ \\ 
$U_{i}^{\prime \prime }=\gamma _{1}\delta _{i-1}+\gamma _{2}\delta
_{i}+\gamma _{1}\delta _{i+1}$%
\end{tabular}
\label{r5}
\end{equation}%
\begin{equation*}
\begin{array}{ll}
\alpha _{1}=\sin ^{2}(\frac{h}{2})\csc (h)\csc (\frac{3h}{2}) & \alpha _{2}=%
\dfrac{2}{1+2\cos (h)} \\ 
\beta _{1}=-\frac{3}{4}\csc (\frac{3h}{2}) & \beta _{2}=\frac{3}{4}\csc (%
\frac{3h}{2}) \\ 
\gamma _{1}=\dfrac{3((1+3\cos (h))\csc ^{2}(\frac{h}{2}))}{16(2\cos (\frac{h%
}{2})+\cos (\frac{3h}{2}))} & \gamma _{2}=-\dfrac{3\cot ^{2}(\frac{h}{2})}{%
2+4\cos (h)}%
\end{array}%
\end{equation*}

The time derivative \ and space derivatives can be approximated by using the
standard finite difference formula and the Crank--Nicolson \ scheme
respectively to have the time-integrated Burger's equation: 
\begin{equation}
\frac{U^{n+1}-U^{n}}{\Delta t}+\frac{(UU_{x})^{n+1}+(UU_{x})^{n}}{2}-\lambda 
\frac{U_{xx}^{n+1}+U_{xx}^{n}}{2}=0  \label{r6}
\end{equation}%
where $U^{n+1}=U(x,t)$ is the solution of the equation at the $(n+1)$th time
level. Here $t^{n+1}$ $=t^{n}+t$, and $\Delta t$ is the time step,
superscripts denote $n$ th time level , $t^{n}=n\Delta t$

The nonlinear term $(UU_{x})^{n+1}$ in Eq. (\ref{r6}) is linearized by using
the following form \cite{rubin,ru1}:%
\begin{equation}
(UUx)^{n+1}=U^{n+1}U_{x}^{n}+U^{n}U_{x}^{n+1}-U^{n}U_{x}^{n}  \label{r7}
\end{equation}%
So linearized time-integrated Burger' equation have the following form:%
\begin{equation}
U^{n+1}-U^{n}+\frac{\Delta t}{2}(U^{n+1}U_{x}^{n}+U^{n}U_{x}^{n+1})-\lambda 
\frac{\Delta t}{2}(U_{xx}^{n+1}-U_{xx}^{n})=0  \label{r8}
\end{equation}%
Substitution \ref{r4} into \ref{r8} \ and evaluation resulting equation at
knots leads to the fully-discretized equation:

\begin{equation}
\begin{tabular}{l}
$\left( \alpha _{1}+\dfrac{\Delta t}{2}\left( \alpha _{1}L_{2}+\beta
_{1}L_{1}-\lambda \gamma _{1}\right) \right) \delta _{m-1}^{n+1}+\left(
\alpha _{2}+\dfrac{\Delta t}{2}\left( \alpha _{2}L_{2}-\lambda \gamma
_{2}\right) \right) \delta _{m}^{n+1}+$ \\ 
$\left( \alpha _{3}+\dfrac{\Delta t}{2}\left( \alpha _{3}L_{2}+\beta
_{2}L_{1}-\lambda \gamma _{3}\right) \right) \delta _{m+1}^{n+1}=(\alpha
_{1}-\lambda \dfrac{\Delta t}{2}\gamma _{1})\delta _{m-1}^{n}+$ \\ 
$(\alpha _{2}-\lambda \dfrac{\Delta t}{2}\gamma _{2})\delta _{m}^{n}+(\alpha
_{3}-\lambda \dfrac{\Delta t}{2}\gamma _{3})\delta _{m+1}^{n}$%
\end{tabular}
\label{r9}
\end{equation}%
where%
\begin{eqnarray*}
L_{1} &=&\alpha _{1}\delta _{i-1}+\alpha _{2}\delta _{i}+\alpha _{3}\delta
_{i+1} \\
L_{2} &=&\beta _{1}\delta _{i-1}+\beta _{2}\delta _{i+1}
\end{eqnarray*}

The system consist of $N+1$ linear equation in $N+3$ unknown parameters $%
\mathbf{d}^{n+1}=(\delta _{-1}^{n+1},\delta _{0}^{n+1},\ldots ,\delta
_{N+1}^{n+1})$. The above system can be made solvable by elimination the
time parameters $\delta _{-1},\delta _{N+1}$with help of the \bigskip
boundary conditions $U(x,a)=U_{0},U(x,b)=U_{N}$ when written as 
\begin{eqnarray}
\delta _{-1} &=&\frac{1}{\alpha _{1}}\left( U_{0}-\alpha _{2}\delta
_{0}-\alpha _{3}\delta _{1}\right) ,  \label{r10} \\
\delta _{N+1} &=&\frac{1}{\alpha _{3}}\left( U_{n}-\alpha _{1}\delta
_{N-1}-\alpha _{2}\delta _{N}\right) .  \notag
\end{eqnarray}%
A variant of Thomas algorithm is used to solve the system.

Initial parameters $d^{0}=\delta _{-1}^{0},\delta _{0}^{0},\ldots ,\delta
_{N+1}^{0}$ must be found to start the iteration process. To do so, initial
condition and boundary values of derivative of initial conditions gives the
following equation

\begin{enumerate}
\item $U_{N}(x_{i},0)$ $=U(x_{i},0),$ $i=0,...,N$

\item $(U_{x})_{N}(x_{0},0)=U^{\prime }(x_{0})$

\item $(U_{x})_{N}(x_{N},0)=U^{\prime }(x_{N}).$
\end{enumerate}

the above system yields an ($N+3$)$\times $($N+3$) matrix system, which can
be solved by use of the Thomas algorithm.

Once the initial parameters $d^{0}$ has been obtained from the initial and
boundary conditions, the recurrence relation gives time evolution of vector $%
d^{n},$from the time evolution of the approximate solution $U_{N}(x,t)$ can
be computed via the equation (\ref{r6}).

\section{Numerical tests}

\textbf{(a)} Analytical solution of the Burger's equation with the problem
sine wave initial condition $U(x,0)=sin(\pi x)$ and boundary conditions $%
U(0,t)=U(1,t)=0$ can be expressed as an infinite series \cite{cole}

\begin{equation}
U(x,t)=\frac{4\pi \lambda \dsum \limits_{j=1}^{\infty }j\mathbf{I}_{j}(\frac{%
1}{2\pi \lambda })\sin (j\pi x)\exp (-j^{2}\pi ^{2}\lambda t)}{\mathbf{I}%
_{0}(\frac{1}{2\pi \lambda })+2\dsum \limits_{j=1}^{\infty }\mathbf{I}_{j}(%
\frac{1}{2\pi \lambda })\cos (j\pi x)\exp (-j^{2}\pi ^{2}\lambda t)}
\label{r11}
\end{equation}%
where $\mathbf{I}_{j}$ are the modified Bessel functions. This problem gives
the decay of sinusoidal disturbance. The convergence of the solution \cite%
{r11} is slow for small values of $\lambda $so that the numerical solutions
of the Burger's equation are looked for. Using the parameters $N=40,$ $%
\Delta t=0.0001,\lambda =1,0.1,0.01,0.00,$ graphical solutions at different
times are depicted in the Figs 1-4 The amplitude of the solution decays as
time pass, seen in Fig 1-2 clearly and \ the sharpness through the right
boundary develops when the smaller viscosities are used. he same incidents
also exist for studies given in the paper \cite{QBSpCol,QuinCol}

\begin{equation*}
\begin{array}{cc}
\begin{tabular}{l}
\FRAME{itbpF}{2.7095in}{2.2537in}{0in}{}{}{fig1.bmp}{\special{language
"Scientific Word";type "GRAPHIC";maintain-aspect-ratio TRUE;display
"USEDEF";valid_file "F";width 2.7095in;height 2.2537in;depth
0in;original-width 2.6671in;original-height 2.2139in;cropleft "0";croptop
"1";cropright "1";cropbottom "0";filename 'Fig1.bmp';file-properties
"XNPEU";}} \\ 
Fig. 1: Solutions at different times \\ 
for $\lambda =1,$ $N=40,$ $\Delta t=0.0001.$%
\end{tabular}
& 
\begin{tabular}{l}
\FRAME{itbpF}{2.7224in}{2.2468in}{0in}{}{}{fig2.bmp}{\special{language
"Scientific Word";type "GRAPHIC";maintain-aspect-ratio TRUE;display
"USEDEF";valid_file "F";width 2.7224in;height 2.2468in;depth
0in;original-width 2.6801in;original-height 2.207in;cropleft "0";croptop
"1";cropright "1";cropbottom "0";filename 'Fig2.bmp';file-properties
"XNPEU";}} \\ 
Fig. 2: Solutions at different times \\ 
for $\lambda =0.1,$ $N=40,$ $\Delta t=0.0001.$%
\end{tabular}%
\end{array}%
\end{equation*}%
\begin{equation*}
\begin{array}{cc}
\begin{tabular}{l}
\FRAME{itbpF}{2.7294in}{2.2935in}{0in}{}{}{fig3.bmp}{\special{language
"Scientific Word";type "GRAPHIC";maintain-aspect-ratio TRUE;display
"USEDEF";valid_file "F";width 2.7294in;height 2.2935in;depth
0in;original-width 2.687in;original-height 2.2537in;cropleft "0";croptop
"1";cropright "1";cropbottom "0";filename 'Fig3.bmp';file-properties
"XNPEU";}} \\ 
Fig. 3: Solutions at different times \\ 
for $\lambda =0.01,$ $N=40,$ $\Delta t=0.0001.$%
\end{tabular}
& 
\begin{tabular}{l}
\FRAME{itbpF}{2.7025in}{2.3203in}{0in}{}{}{fig4.bmp}{\special{language
"Scientific Word";type "GRAPHIC";maintain-aspect-ratio TRUE;display
"USEDEF";valid_file "F";width 2.7025in;height 2.3203in;depth
0in;original-width 2.6593in;original-height 2.2805in;cropleft "0";croptop
"1";cropright "1";cropbottom "0";filename 'Fig4.bmp';file-properties
"XNPEU";}} \\ 
Fig. 4: Solutions at different times \\ 
for $\lambda =0.001,$ $N=40,$ $\Delta t=0.0001.$%
\end{tabular}%
\end{array}%
\end{equation*}%
The results of proposed numerical methods are compared with \ the cubic
B-spline collocation, cubic B-spline Galerkin. Galerkin procedure are seen
to produce slightly same results with the CTB collocation method. Our
advantage is that the cost of the CTB procedure is less than \ the Galerkin
methods given in the tables 1-3. 
\begin{equation*}
\begin{tabular}{|l|}
\hline
$\text{Table 2: Comparison of the numerical solutions of Problem 1}$ \\ 
\hline
$\text{obtained for }\lambda =1.$and $N=40,$ $\Delta t=0.0001$ \\ \hline
at different times with the exact solutions \\ \hline
$%
\begin{tabular}{llllll}
$x$ & $t$ & $%
\begin{array}{c}
\text{Present} \\ 
\text{ }%
\end{array}%
$ & $%
\begin{array}{c}
\text{Ref.\cite{CBSpCol}} \\ 
(N=80)%
\end{array}%
$ & $%
\begin{array}{c}
\text{\cite{BSpGal}} \\ 
\text{ }%
\end{array}%
$ & $%
\begin{array}{c}
\text{Exact} \\ 
\text{ }%
\end{array}%
$ \\ \hline
0.25 & 0.4 & \multicolumn{1}{c}{$0.01355$} & \multicolumn{1}{c}{$0.01357$} & 
\multicolumn{1}{c}{$0.01357$} & \multicolumn{1}{c}{$0.01357$} \\ 
& 0.6 & \multicolumn{1}{c}{$0.00188$} & \multicolumn{1}{c}{$0.00189$} & 
\multicolumn{1}{c}{$0.00189$} & \multicolumn{1}{c}{$0.00189$} \\ 
& 0.8 & \multicolumn{1}{c}{$0.00026$} & \multicolumn{1}{c}{$0.00026$} & 
\multicolumn{1}{c}{$0.00026$} & \multicolumn{1}{c}{$0.00026$} \\ 
& 1.0 & \multicolumn{1}{c}{$0.00004$} & \multicolumn{1}{c}{$0.00004$} & 
\multicolumn{1}{c}{$0.00004$} & \multicolumn{1}{c}{$0.00004$} \\ 
& 3.0 & \multicolumn{1}{c}{$0.00000$} & \multicolumn{1}{c}{$0.00000$} & 
\multicolumn{1}{c}{$0.00000$} & \multicolumn{1}{c}{$0.00000$} \\ 
&  & \multicolumn{1}{c}{} & \multicolumn{1}{c}{} & \multicolumn{1}{c}{} & 
\multicolumn{1}{c}{} \\ 
0.50 & 0.4 & \multicolumn{1}{c}{$0.01920$} & \multicolumn{1}{c}{$0.01923$} & 
\multicolumn{1}{c}{$0.01924$} & \multicolumn{1}{c}{$0.01924$} \\ 
& 0.6 & \multicolumn{1}{c}{$0.00266$} & \multicolumn{1}{c}{$0.00267$} & 
\multicolumn{1}{c}{$0.00267$} & \multicolumn{1}{c}{$0.00267$} \\ 
& 0.8 & \multicolumn{1}{c}{$0.00037$} & \multicolumn{1}{c}{$0.00037$} & 
\multicolumn{1}{c}{$0.00037$} & \multicolumn{1}{c}{$0.00037$} \\ 
& 1.0 & \multicolumn{1}{c}{$0.00005$} & \multicolumn{1}{c}{$0.00005$} & 
\multicolumn{1}{c}{$0.00005$} & \multicolumn{1}{c}{$0.00005$} \\ 
& 3.0 & \multicolumn{1}{c}{$0.00000$} & \multicolumn{1}{c}{$0.00000$} & 
\multicolumn{1}{c}{$0.00000$} & \multicolumn{1}{c}{$0.00000$} \\ 
&  & \multicolumn{1}{c}{} & \multicolumn{1}{c}{} & \multicolumn{1}{c}{} & 
\multicolumn{1}{c}{} \\ 
0.75 & 0.4 & \multicolumn{1}{c}{$0.01361$} & \multicolumn{1}{c}{$0.01362$} & 
\multicolumn{1}{c}{$0.01363$} & \multicolumn{1}{c}{$0.01363$} \\ 
& 0.6 & \multicolumn{1}{c}{$0.00188$} & \multicolumn{1}{c}{$0.00189$} & 
\multicolumn{1}{c}{$0.00189$} & \multicolumn{1}{c}{$0.00189$} \\ 
& 0.8 & \multicolumn{1}{c}{$0.00026$} & \multicolumn{1}{c}{$0.00026$} & 
\multicolumn{1}{c}{$0.00026$} & \multicolumn{1}{c}{$0.00026$} \\ 
& 1.0 & \multicolumn{1}{c}{$0.00004$} & \multicolumn{1}{c}{$0.00004$} & 
\multicolumn{1}{c}{$0.00004$} & \multicolumn{1}{c}{$0.00004$} \\ 
& 3.0 & \multicolumn{1}{c}{$0.00000$} & \multicolumn{1}{c}{$0.00000$} & 
\multicolumn{1}{c}{$0.00000$} & \multicolumn{1}{c}{$0.00000$}%
\end{tabular}%
$ \\ \hline
\end{tabular}%
\end{equation*}%
\begin{equation*}
\begin{tabular}{|l|}
\hline
$\text{Table 3: Comparison of the numerical solutions of Problem 1}$ \\ 
\hline
$\text{obtained for }\lambda =0.1.$and $N=40,$ $\Delta t=0.0001$ \\ \hline
at different times with the exact solutions \\ \hline
\multicolumn{1}{|c|}{$%
\begin{tabular}{lllllll}
$x$ & $t$ & $%
\begin{array}{c}
\text{Present} \\ 
\text{ }%
\end{array}%
$ & $%
\begin{array}{c}
\text{Ref.\cite{CBSpCol}} \\ 
(N=80)%
\end{array}%
$ & $%
\begin{array}{c}
\text{Ref.\cite{QBSpCol2}} \\ 
\text{ }%
\end{array}%
$ & $%
\begin{array}{c}
\text{Ref.\cite{BSpGal}} \\ 
\text{ }%
\end{array}%
$ & $%
\begin{array}{c}
\text{Exact} \\ 
\text{ }%
\end{array}%
$ \\ \hline
0.25 & 0.4 & \multicolumn{1}{c}{$0.30892$} & \multicolumn{1}{c}{$0.30890$} & 
\multicolumn{1}{c}{$0.30891$} & \multicolumn{1}{c}{$0.30890$} & 
\multicolumn{1}{c}{$0.30889$} \\ 
& 0.6 & \multicolumn{1}{c}{$0.24078$} & \multicolumn{1}{c}{$0.24075$} & 
\multicolumn{1}{c}{$0.24075$} & \multicolumn{1}{c}{$0.24074$} & 
\multicolumn{1}{c}{$0.24074$} \\ 
& 0.8 & \multicolumn{1}{c}{$0.19572$} & \multicolumn{1}{c}{$0.19569$} & 
\multicolumn{1}{c}{$0.19568$} & \multicolumn{1}{c}{$0.19568$} & 
\multicolumn{1}{c}{$0.19568$} \\ 
& 1.0 & \multicolumn{1}{c}{$0.16261$} & \multicolumn{1}{c}{$0.16258$} & 
\multicolumn{1}{c}{$0.16257$} & \multicolumn{1}{c}{$0.16257$} & 
\multicolumn{1}{c}{$0.16256$} \\ 
& 3.0 & \multicolumn{1}{c}{$0.02718$} & \multicolumn{1}{c}{$0.02720$} & 
\multicolumn{1}{c}{$0.02721$} & \multicolumn{1}{c}{$0.02720$} & 
\multicolumn{1}{c}{$0.02720$} \\ 
&  & \multicolumn{1}{c}{} & \multicolumn{1}{c}{} & \multicolumn{1}{c}{} & 
\multicolumn{1}{c}{} & \multicolumn{1}{c}{} \\ 
0.50 & 0.4 & \multicolumn{1}{c}{$0.56971$} & \multicolumn{1}{c}{$0.56965$} & 
\multicolumn{1}{c}{$0.56969$} & \multicolumn{1}{c}{$0.56964$} & 
\multicolumn{1}{c}{$0.56963$} \\ 
& 0.6 & \multicolumn{1}{c}{$0.44730$} & \multicolumn{1}{c}{$0.44723$} & 
\multicolumn{1}{c}{$0.44723$} & \multicolumn{1}{c}{$0.44721$} & 
\multicolumn{1}{c}{$0.44721$} \\ 
& 0.8 & \multicolumn{1}{c}{$0.35932$} & \multicolumn{1}{c}{$0.35925$} & 
\multicolumn{1}{c}{$0.35926$} & \multicolumn{1}{c}{$0.35924$} & 
\multicolumn{1}{c}{$0.35924$} \\ 
& 1.0 & \multicolumn{1}{c}{$0.29197$} & \multicolumn{1}{c}{$0.29192$} & 
\multicolumn{1}{c}{$0.29193$} & \multicolumn{1}{c}{$0.29191$} & 
\multicolumn{1}{c}{$0.29192$} \\ 
& 3.0 & \multicolumn{1}{c}{$0.04017$} & \multicolumn{1}{c}{$0.04019$} & 
\multicolumn{1}{c}{$0.04021$} & \multicolumn{1}{c}{$0.04020$} & 
\multicolumn{1}{c}{$0.04021$} \\ 
&  & \multicolumn{1}{c}{} & \multicolumn{1}{c}{} & \multicolumn{1}{c}{} & 
\multicolumn{1}{c}{} & \multicolumn{1}{c}{} \\ 
0.75 & 0.4 & \multicolumn{1}{c}{$0.62524$} & \multicolumn{1}{c}{$0.62538$} & 
\multicolumn{1}{c}{$0.62543$} & \multicolumn{1}{c}{$0.62541$} & 
\multicolumn{1}{c}{$0.62544$} \\ 
& 0.6 & \multicolumn{1}{c}{$0.48698$} & \multicolumn{1}{c}{$0.48715$} & 
\multicolumn{1}{c}{$0.48723$} & \multicolumn{1}{c}{$0.48719$} & 
\multicolumn{1}{c}{$0.48721$} \\ 
& 0.8 & \multicolumn{1}{c}{$0.37369$} & \multicolumn{1}{c}{$0.37385$} & 
\multicolumn{1}{c}{$0.37394$} & \multicolumn{1}{c}{$0.37390$} & 
\multicolumn{1}{c}{$0.37392$} \\ 
& 1.0 & \multicolumn{1}{c}{$0.28727$} & \multicolumn{1}{c}{$0.28741$} & 
\multicolumn{1}{c}{$0.28750$} & \multicolumn{1}{c}{$0.28746$} & 
\multicolumn{1}{c}{$0.28747$} \\ 
& 3.0 & \multicolumn{1}{c}{$0.02974$} & \multicolumn{1}{c}{$0.02976$} & 
\multicolumn{1}{c}{$0.02978$} & \multicolumn{1}{c}{$0.02977$} & 
\multicolumn{1}{c}{$0.02977$}%
\end{tabular}%
$} \\ \hline
\end{tabular}%
\end{equation*}%
\begin{equation*}
\begin{tabular}{|l|}
\hline
$\text{Table 4: Comparison of the numerical solutions of Problem 1}$ \\ 
\hline
$\text{obtained for }\lambda =0.01.$and $N=40,$ $\Delta t=0.0001$ \\ \hline
at different times with the exact solutions \\ \hline
$%
\begin{tabular}{lllllll}
$x$ & $t$ & $%
\begin{array}{c}
\text{Present} \\ 
\text{ }%
\end{array}%
$ & $%
\begin{array}{c}
\text{Ref.\cite{CBSpCol}} \\ 
(N=80)%
\end{array}%
$ & $%
\begin{array}{c}
\text{Ref.\cite{QBSpCol2}} \\ 
\text{ }%
\end{array}%
$ & $%
\begin{array}{c}
\text{Ref.\cite{BSpGal}} \\ 
\text{ }%
\end{array}%
$ & $%
\begin{array}{c}
\text{Exact} \\ 
\text{ }%
\end{array}%
$ \\ \hline
0.25 & 0.4 & \multicolumn{1}{c}{$0.34191$} & \multicolumn{1}{c}{$0.34192$} & 
\multicolumn{1}{c}{$0.34192$} & \multicolumn{1}{c}{$0.34192$} & 
\multicolumn{1}{c}{$0.34191$} \\ 
& 0.6 & \multicolumn{1}{c}{$0.26896$} & \multicolumn{1}{c}{$0.26897$} & 
\multicolumn{1}{c}{$0.22894$} & \multicolumn{1}{c}{$0.26897$} & 
\multicolumn{1}{c}{$0.22896$} \\ 
& 0.8 & \multicolumn{1}{c}{$0.22148$} & \multicolumn{1}{c}{$0.22148$} & 
\multicolumn{1}{c}{$0.22144$} & \multicolumn{1}{c}{$0.22148$} & 
\multicolumn{1}{c}{$0.22148$} \\ 
& 1.0 & \multicolumn{1}{c}{$0.18819$} & \multicolumn{1}{c}{$0.18819$} & 
\multicolumn{1}{c}{$0.18815$} & \multicolumn{1}{c}{$0.18819$} & 
\multicolumn{1}{c}{$0.18819$} \\ 
& 3.0 & \multicolumn{1}{c}{$0.07511$} & \multicolumn{1}{c}{$0.07511$} & 
\multicolumn{1}{c}{$0.07509$} & \multicolumn{1}{c}{$0.07511$} & 
\multicolumn{1}{c}{$0.07511$} \\ 
&  & \multicolumn{1}{c}{} & \multicolumn{1}{c}{} & \multicolumn{1}{c}{} & 
\multicolumn{1}{c}{} & \multicolumn{1}{c}{} \\ 
0.50 & 0.4 & \multicolumn{1}{c}{$0.66071$} & \multicolumn{1}{c}{$0.66071$} & 
\multicolumn{1}{c}{$0.66075$} & \multicolumn{1}{c}{$0.66071$} & 
\multicolumn{1}{c}{$0.66071$} \\ 
& 0.6 & \multicolumn{1}{c}{$0.52942$} & \multicolumn{1}{c}{$0.52942$} & 
\multicolumn{1}{c}{$0.52946$} & \multicolumn{1}{c}{$0.52942$} & 
\multicolumn{1}{c}{$0.52942$} \\ 
& 0.8 & \multicolumn{1}{c}{$0.43914$} & \multicolumn{1}{c}{$0.43914$} & 
\multicolumn{1}{c}{$0.43917$} & \multicolumn{1}{c}{$0.43914$} & 
\multicolumn{1}{c}{$0.43914$} \\ 
& 1.0 & \multicolumn{1}{c}{$0.37442$} & \multicolumn{1}{c}{$0.37442$} & 
\multicolumn{1}{c}{$0.37444$} & \multicolumn{1}{c}{$0.37442$} & 
\multicolumn{1}{c}{$0.37442$} \\ 
& 3.0 & \multicolumn{1}{c}{$0.15017$} & \multicolumn{1}{c}{$0.15018$} & 
\multicolumn{1}{c}{$0.15016$} & \multicolumn{1}{c}{$0.15018$} & 
\multicolumn{1}{c}{$0.15018$} \\ 
&  & \multicolumn{1}{c}{} & \multicolumn{1}{c}{} & \multicolumn{1}{c}{} & 
\multicolumn{1}{c}{} & \multicolumn{1}{c}{} \\ 
0.75 & 0.4 & \multicolumn{1}{c}{$0.91029$} & \multicolumn{1}{c}{$0.91027$} & 
\multicolumn{1}{c}{$0.91023$} & \multicolumn{1}{c}{$0.91027$} & 
\multicolumn{1}{c}{$0.91026$} \\ 
& 0.6 & \multicolumn{1}{c}{$0.76725$} & \multicolumn{1}{c}{$0.76725$} & 
\multicolumn{1}{c}{$0.76728$} & \multicolumn{1}{c}{$0.76724$} & 
\multicolumn{1}{c}{$0.76724$} \\ 
& 0.8 & \multicolumn{1}{c}{$0.64740$} & \multicolumn{1}{c}{$0.64740$} & 
\multicolumn{1}{c}{$0.64744$} & \multicolumn{1}{c}{$0.64740$} & 
\multicolumn{1}{c}{$0.64740$} \\ 
& 1.0 & \multicolumn{1}{c}{$0.55605$} & \multicolumn{1}{c}{$0.55605$} & 
\multicolumn{1}{c}{$0.55609$} & \multicolumn{1}{c}{$0.55605$} & 
\multicolumn{1}{c}{$0.55605$} \\ 
& 3.0 & \multicolumn{1}{c}{$0.22489$} & \multicolumn{1}{c}{$0.22483$} & 
\multicolumn{1}{c}{$0.22481$} & \multicolumn{1}{c}{$0.22481$} & 
\multicolumn{1}{c}{$0.22481$}%
\end{tabular}%
$ \\ \hline
\end{tabular}%
\end{equation*}

\textbf{(b) }Well-known other solution of the Burger's equation is%
\begin{equation}
U(x,t)=\dfrac{\alpha +\mu +(\mu -\alpha )\exp \eta }{1+\exp \eta },\text{ }%
0\leq x\leq 1,\text{ }t\geq 0,  \label{r12}
\end{equation}%
where $\eta =\dfrac{\alpha (x-\mu t-\gamma )}{\lambda }.$ $\alpha ,~\mu $
and $\gamma $ are arbitrary constants. Initial conditions are$U(0,t)=1,$ $%
U(1,t)=0.2$ or $U_{x}(0,t)=0,$ $U_{x}(1,t)=0,$ for $t\geq 0$ \ This form of
the solution is known as the travelling wave equation and respresent the
propogation of the wave front through the right. Parameter $\lambda $
determine the sharpness of the solution

The initila solutions are taken from the analytical solution when $t=0$. The
program is run for the parameters $\alpha =0.4,$ $\mu =0.6,$ $\gamma =0.125$
and $\lambda =0.01,h=1/36,\Delta t=0.001$. solutions at some space values $x$
are presented in Table 5 and compared with those obtained in the studies 
\cite{CBSpCol,BSpGal,BSpGal} using Cubic B-spline collocation,
quadratic/Cubic B-spline Galerkin methods. Solution behaviours are
illustrated in Fig 5-6 for the coefficient $\lambda =0.01$ and $0.001$ at
times $t=0,0.4,0.8,1.2.$ With smaller $\lambda =0.001,$the sharp front is
formed and propogates to right during run of the program. Graphical
presentation of the absolute errors at time $t=0.4$ is drawn im Figs 7-8$\ $

\begin{equation*}
\begin{tabular}{|l|}
\hline
Table 5: Comparison of results at time \\ \hline
$t=0.5,$ $h=1/36,$ $\Delta t=0.01,$ $\lambda =0.01$ \\ \hline
$%
\begin{array}{lccccc}
x & 
\begin{array}{c}
\text{Present} \\ 
\text{~}%
\end{array}
& 
\begin{array}{c}
\text{Ref. \cite{CBSpCol}} \\ 
\Delta t=0.025%
\end{array}
& 
\begin{array}{c}
\text{Ref. \cite{BSpGal}} \\ 
\text{(QBGM)}%
\end{array}
& 
\begin{array}{c}
\text{Ref. \cite{BSpGal}} \\ 
\text{(CBGM)}%
\end{array}
& 
\begin{array}{c}
\text{Exact} \\ 
\text{ }%
\end{array}
\\ 
{\small 0.000} & {\small 1.} & {\small 1.} & {\small 1.} & {\small 1.} & 
{\small 1.} \\ 
{\small 0.056} & {\small 1.} & {\small 1.} & {\small 1.} & {\small 1.} & 
{\small 1.} \\ 
{\small 0.111} & {\small 1.} & {\small 1.} & {\small 1.} & {\small 1.} & 
{\small 1.} \\ 
{\small 0.167} & {\small 1.} & {\small 1.} & {\small 1.} & {\small 1.} & 
{\small 1.} \\ 
{\small 0.222} & {\small 1.} & {\small 1.} & {\small 1.} & {\small 1.} & 
{\small 1.} \\ 
{\small 0.278} & {\small 0.999} & {\small 0.999} & {\small 0.998} & {\small %
0.998} & {\small 0.998} \\ 
{\small 0.333} & {\small 0.983} & {\small 0.986} & {\small 0.980} & {\small %
0.980} & {\small 0.980} \\ 
{\small 0.389} & {\small 0.845} & {\small 0.850} & {\small 0.841} & {\small %
0.842} & {\small 0.847} \\ 
{\small 0.444} & {\small 0.456} & {\small 0.448} & {\small 0.458} & {\small %
0.457} & {\small 0.452} \\ 
{\small 0.500} & {\small 0.237} & {\small 0.236} & {\small 0.240} & {\small %
0.241} & {\small 0.238} \\ 
{\small 0.556} & {\small 0.203} & {\small 0.204} & {\small 0.205} & {\small %
0.205} & {\small 0.204} \\ 
{\small 0.611} & {\small 0.2} & {\small 0.2} & {\small 0.2} & {\small 0.2} & 
{\small 0.2} \\ 
{\small 0.667} & {\small 0.2} & {\small 0.2} & {\small 0.2} & {\small 0.2} & 
{\small 0.2} \\ 
{\small 0.722} & {\small 0.2} & {\small 0.2} & {\small 0.2} & {\small 0.2} & 
{\small 0.2} \\ 
{\small 0.778} & {\small 0.2} & {\small 0.2} & {\small 0.2} & {\small 0.2} & 
{\small 0.2} \\ 
{\small 0.833} & {\small 0.2} & {\small 0.2} & {\small 0.2} & {\small 0.2} & 
{\small 0.2} \\ 
{\small 0.889} & {\small 0.2} & {\small 0.2} & {\small 0.2} & {\small 0.2} & 
{\small 0.2} \\ 
{\small 0.944} & {\small 0.2} & {\small 0.2} & {\small 0.2} & {\small 0.2} & 
{\small 0.2} \\ 
{\small 1.000} & {\small 0.2} & {\small 0.2} & {\small 0.2} & {\small 0.2} & 
{\small 0.2}%
\end{array}%
$ \\ \hline
\end{tabular}%
\end{equation*}%
\begin{equation*}
\begin{array}{c}
\FRAME{itbpF}{2.7224in}{2.2537in}{0in}{}{}{fig5.bmp}{\special{language
"Scientific Word";type "GRAPHIC";maintain-aspect-ratio TRUE;display
"USEDEF";valid_file "F";width 2.7224in;height 2.2537in;depth
0in;original-width 2.6801in;original-height 2.2139in;cropleft "0";croptop
"1";cropright "1";cropbottom "0";filename 'Fig5.bmp';file-properties
"XNPEU";}} \\ 
\multicolumn{1}{l}{\text{Figure 5: Solutions at different times for}} \\ 
\multicolumn{1}{l}{\lambda =0.01,\text{ }h=1/36,\text{ }\Delta t=0.001,\text{
}x\in \lbrack 0,1].}%
\end{array}%
\begin{tabular}{l}
\FRAME{itbpF}{2.7224in}{2.3065in}{0in}{}{}{fig6.bmp}{\special{language
"Scientific Word";type "GRAPHIC";maintain-aspect-ratio TRUE;display
"USEDEF";valid_file "F";width 2.7224in;height 2.3065in;depth
0in;original-width 2.6801in;original-height 2.2667in;cropleft "0";croptop
"1";cropright "1";cropbottom "0";filename 'Fig6.bmp';file-properties
"XNPEU";}} \\ 
Figure 6: Solutions at different times for \\ 
$\lambda =0.005,$ $h=1/36,$ $\Delta t=0.001,\text{ }x\in \lbrack 0,1].$%
\end{tabular}%
\end{equation*}%
\begin{equation*}
\begin{array}{cc}
\begin{tabular}{l}
\FRAME{itbpF}{2.789in}{2.2399in}{0in}{}{}{fig7.bmp}{\special{language
"Scientific Word";type "GRAPHIC";maintain-aspect-ratio TRUE;display
"USEDEF";valid_file "F";width 2.789in;height 2.2399in;depth
0in;original-width 2.7466in;original-height 2.2001in;cropleft "0";croptop
"1";cropright "1";cropbottom "0";filename 'Fig7.bmp';file-properties
"XNPEU";}} \\ 
$\text{Fig. 7}$: $L_{\infty }$ error norm for $\lambda =0.01$ and $h=1/36$%
\end{tabular}
& 
\begin{tabular}{l}
\FRAME{itbpF}{2.783in}{2.2329in}{0in}{}{}{fig8.bmp}{\special{language
"Scientific Word";type "GRAPHIC";maintain-aspect-ratio TRUE;display
"USEDEF";valid_file "F";width 2.783in;height 2.2329in;depth
0in;original-width 2.7397in;original-height 2.1932in;cropleft "0";croptop
"1";cropright "1";cropbottom "0";filename 'Fig8.bmp';file-properties
"XNPEU";}} \\ 
$\text{Fig. 8 }L_{\infty }$ error norm $\lambda =0.005$ and $h=1/36$%
\end{tabular}%
\end{array}%
\end{equation*}%
The collocation methods with trigonometric B-spline functions is made up to
find solutions the Burger's equation. We have hown that methods is capable
of producing solutions of the Burgers equation fairly. The method can be
used as an alternative to the methods accompanied B-spline functions.

\end{document}